\magnification1200

\font\frak=eufm10
\def\g{\hbox{\frak\char'0147}}

\def\s{\hbox{\frak\char'0163}}
\def\p{\hbox{\frak\char'0160}}

\def\l{\hbox{\frak\char'0154}}

\def\h{\hbox{\frak\char'0150}}

\def\n{\hbox{\frak\char'0156}}
\def\1{\hbox{\frak\char'0061}}

\font\block=msbm10
\def\C{\hbox{\block\char'0103}}

\def\Z{\hbox{\block\char'0132}}

\def\+{\hbox{\block\char'0156}}
\def\H{\hbox{\block\char'0110}}
\def\-{\hbox{\block\char'0157}}

\font\red=eusm10
\def\F{\hbox{\red\char'0106}}
\def\S{\hbox{\red\char'0123}}
\def\L{\hbox{\red\char'0114}}
\font\green=msam10
\def\SQ{\hbox{\green\char'0003}}

\font\white=cmbsy10
\def\H1{\hbox{\white\char'0110}}
\def\E1{\hbox{\white\char'0105}}
\def\F1{\hbox{\white\char'0106}}

\baselineskip=1.3\normalbaselineskip
\font\bigtenrm=cmr10 scaled\magstep1

\font\tenbf=cmbx10 scaled\magstep1
\font\bf=cmbx10 scaled\magstep0

\font\new=eurm10
\def\x{\hbox{\new\char'0170}}
\def\y{\hbox{\new\char'0171}}
\def\p1{\hbox{\new\char'0160}}
\def\h1{\hbox{\new\char'0150}}
\def\centre{\hbox{\new\char'0103}}

\font\frak=eufm10

\font\ninerm=cmr9

\centerline{\bf  SEMI-INFINITE COHOMOLOGY AND SUPERCONFORMAL ALGEBRAS}

\vskip 0.2in
\centerline{\bigtenrm By Elena Poletaeva}
\vskip 0.2in
{\ninerm  {\sl Abstract} - 
We describe representations of
certain  superconformal algebras in the
semi-infinite Weil complex related to the loop algebra
of a complex finite-dimensional Lie algebra and in the 
semi-infinite cohomology.
We show that in the case where the Lie algebra
is endowed with a non-degenerate invariant symmetric bilinear form,
the relative semi-infinite cohomology
of the loop algebra has a structure, 
which is analogous to the classical structure of the de Rham
cohomology in K\"ahler geometry.}
\vskip 0.2in
\centerline{\bf  COHOMOLOGIE SEMI-INFINIE ET ALG\`EBRES SUPERCONFORMES}

\vskip 0.2in
{\ninerm {\sl R\'esum\'e} -  
Nous d\'ecrivons des repr\'esentations
de certaines alg\`ebres superconformes dans le
complexe de Weil semi-infini de l'alg\`ebre des lacets
d'une alg\`ebre de Lie complexe de dimension finie et
dans la cohomologie semi-infinie.
Nous  d\'emontrons que dans le cas o\`u l'alg\`ebre de Lie
est munie d'une forme sym\'etrique non d\'eg\'en\'er\'ee invariante,
la cohomologie semi-infinie relative
de l'alg\`ebre des lacets admet une structure, qui est
l'analogue de la structure classique
de la cohomologie de de Rham des vari\'et\'es k\"ahl\'eriennes.}

\vskip 0.2in
\centerline{\tenbf 1. Introduction}
\vskip 0.2in

B. Feigin and E. Frenkel have introduced 
a semi-infinite analogue of the
Weil complex based on the space
$$W^{{\infty \over 2} + *}(\g) = 
S^{{\infty \over 2} + *}(\g)
\otimes \Lambda^{{\infty \over 2} + *}(\g).\eqno (1.1)$$
In their construction $\g = \oplus_{n\in\Z}\g_n$
is a graded Lie algebra,
$S^{{\infty \over 2} + *}(\g)$ and 
$\Lambda^{{\infty \over 2} + *}(\g)$ are some semi-infinite
analogues of the symmetric and exterior power modules, [FF].
As in the classical case,
two differentials, $d$ and $\h1$, are defined on 
$W^{{\infty \over 2} + *}(\g)$. They are analogous to the 
differential in Lie algebra (co)homology and 
the Koszul differential, respectively. The semi-infinite Weil complex
$$\lbrace W^{{\infty \over 2} + *}(\g), \hbox{ }d + \h1 \rbrace
 \eqno
(1.2)$$
is acyclic similarly to the classical Weil complex.
The cohomology  of the complex
$$\lbrace W^{{\infty \over 2} + *}(\g), \hbox{ }d  \rbrace
 \eqno
(1.3)$$
is called the 
  {\it semi-infinite cohomology} of $\g$
with coefficients in its ``adjoint semi-infinite symmetric powers"
$H^{{\infty \over 2} + *}(\g, S^{{\infty \over 2} + *}(\g))$.
One can also define  {\it the relative semi-infinite Weil complex}
$W^{{\infty \over 2} + *}_{rel}(\g)$
(relatively $\g_0$), and {\it the  relative semi-infinite cohomology}
\hfil\break
$H^{{\infty \over 2} + *}(\g, \g_0, S^{{\infty \over 2} + *}(\g))$, 
[FF].

E. Getzler has shown that the semi-infinite Weil complex  of the Virasoro
algebra admits an action of the {\it $N = 2$ superconformal algebra}, [G].

Recall that 
a {\it superconformal algebra} (SCA)  is a simple complex
Lie superalgebra $\s$, such that it 
contains the centerless Virasoro algebra (i.e. Witt algebra)
 $Witt = \oplus_{n\in \Z}\C L_n$
as a subalgebra, and has growth 1. The $\Z$-graded 
superconformal algebras are ones for which
$ad L_0$ is diagonalizable with finite-dimensional eigenspaces, [KL]:
$$\s = \oplus_{j}\s_j, \s_j = \lbrace x\in\s \mid [L_0, x] = jx \rbrace.
\eqno(1.4)$$

In this work we consider the semi-infinite Weil complex constructed for
the next natural (after the Virasoro algebra) class of graded Lie algebras:
the loop algebras of the complex finite-dimensional Lie algebras.
The action of the Virasoro algebra on such complex 
is ensured by the fact that it has a structure
of a vertex operator superalgebra (see [Ak]).

Let $\g$
be a complex
finite-dimensional Lie algebra, and
$\tilde{\g} = \g\otimes \C[t, t^{-1}]$ be the corresponding loop algebra.
We obtain a representation of the   $N = 2$ SCA
in the semi-infinite Weil complex
$W^{{\infty \over 2} + *}(\tilde{\g})$ 
 and in the  semi-infinite cohomology
$H^{{\infty \over 2} + *}(\tilde{\g}, S^{{\infty \over
2} + *}(\tilde{\g}))$ with central charge $3dim \g$.
We extend the representation of the  $N = 2$ SCA in 
$W^{{\infty \over 2} + *}(\tilde{\g})$ to a representation of the
one-parameter family $\hat{S}'(2, \alpha)$ of
deformations of the  $N = 4$ SCA (see [Ad] and [KL]).
In the case, when $\g$ is 
endowed with a non-degenerate invariant symmetric
bilinear form, we obtain a representation 
of $\hat S'(2, 0)$  in 
$H^{{\infty \over 2} + *}(\tilde{\g},    S^{{\infty \over 2} + *}(\tilde{\g}))$.
Finally, there exists a representation
of a central extension of
the Lie superalgebra of all derivations of $S'(2, 0)$  
in the relative semi-infinite cohomology  
$H^{{\infty \over 2} + *}(\tilde{\g},
\tilde{\g}_0, S^{{\infty \over 2} + *}(\tilde{\g}))$.

It was shown in [FGZ] that the cohomology of the 
relative semi-infinite complex
$C^*_{\infty}(\l, \l_0, V)$, where $\l$ is a complex graded Lie algebra, and
$V$ is a graded Hermitian $\l$-module, has 
(under certain conditions) a structure
analogous to that of de Rham cohomology
in K\"ahler geometry.

Recall that given a compact K\"ahler
manifold $M$,
there exist 
a number of classical operators on the
space of differential forms on $M$,
such as the differentials
$\partial, \bar{\partial}, d, d_c$, their corresponding adjoint 
operators and the associated Laplacians (see [GH]). 
There also exists an action of $\s\l(2)$ on $H^*(M)$
according to the Lefschetz theorem. All these operators satisfy
a series of identities known as Hodge identities, [GH].
Naturally, the classical operators form a finite-dimensional Lie superalgebra.

We show that given a complex finite-dimensional
Lie algebra $\g$ endowed with 
a non-degenerate invariant symmetric bilinear form,
there exist the analogues of the classical operators on the complex
$W^{{\infty \over 2} + *}_{rel}(\tilde{\g})$. 
We prove that  the exterior derivations of
$S'(2, 0)$ form an $\s\l(2)$, and observe that they define
an $\s\l(2)$-module structure on 
$H^{{\infty \over 2} + *}(\tilde{\g},
\tilde{\g}_0, S^{{\infty \over 2} + *}(\tilde{\g}))$,
which is the analogue of the 
$\s\l(2)$-module  structure on 
the de Rham cohomology in K\"ahler geometry.

The action of $\hat{S}'(2, 0)$
provides $H^{{\infty \over 2} + *}(\tilde{\g},
\tilde{\g}_0, S^{{\infty \over 2} + *}(\tilde{\g}))$
with eight  series of quadratic operators.
In particular, they
include the semi-infinite Koszul differential $\h1$,
and the semi-infinite analogue of the homotopy operator (cf. [Fu]).
We prove that
the  degree zero part
of the $\Z$-grading of $S'(2, 0)$ defined by the element $L_0\in Witt$, 
is isomorphic to
the Lie superalgebra of classical operators
in K\"ahler geometry.

It would be interesting to interpret the superconformal algebra
$S'(2, 0)$ as ``affinization'' of the classical operators in the  case
of an infinite-dimensional manifold.

This work is partly based on [P1]-[P3].

\vskip 0.2in
\centerline{\tenbf 2. Semi-infinite Weil complex}
\vskip 0.2in

The semi-infinite Weil complex of a graded Lie algebra
was introduced by B. Feigin and E. Frenkel
in [FF]. Recall the necessary definitions. More generally,
let $V = \oplus_{n \in \Z}V_n$ be a graded vector space over $\C$,
 such that
$dim V_n < \infty$. Let
$V' = \oplus_{n \in \Z}V_n'$ be the restricted dual of $V$.
The linear space $V \oplus V'$ carries non-degenerate skew-symmetric and
 symmetric bilinear
forms: $(\cdot ,\cdot )$ and $\lbrace\cdot  , \cdot \rbrace $.
Let $H(V)$ and $C(V)$ be the quotients of the tensor algebra
$T^*(V \oplus V')$ by the ideals generated by the elements of the form
$xy - yx - (x,y)$ and $xy + yx - \lbrace x,  y\rbrace$, respectively,
where $x, y \in V \oplus V' $. We fix $K \in \Z$.
Let $V = V_+ \oplus V_-$ be the corresponding polarization of $V$:
$V_+ = \oplus_{n>K}V_n$, $V_- = \oplus_{n\leq K}V_n$.

The symmetric algebra $S^*(V_+ \oplus V_-')$ is a subalgebra
 of $H(V)$ and
the exterior algebra
$\Lambda^*(V_+ \oplus V_-')$ is a subalgebra of $C(V)$.
Let $S^{{\infty \over 2} + *}(V)$, $\Lambda^{{\infty \over 2} + *}(V)$
be the representations of $H(V)$ and $C(V)$ induced from the trivial
representations $<\1_S>$ and  $<\1_{\Lambda}>$ of
\hfil\break
$S^*(V_+ \oplus V_-')$  and of $\Lambda^*(V_+ \oplus V_-')$, respectively.
Thus we obtain some semi-infinite analogues of symmetric and
 exterior power
modules.
Denote the actions of $H(V)$ and $C(V)$ on these modules by
$\beta (x), \gamma (x')$ and $\tau (x), \varepsilon (x')$, respectively,
for $x\in V$, $x' \in V'$. Notice that each element of
$S^{{\infty \over 2} + *}(V)$ and of
 $\Lambda^{{\infty \over 2} + *}(V)$
is a finite linear combination of the monomials of the type
$\gamma (x_1')\ldots \gamma (x_k')\beta (y_1)\ldots
\beta (y_m)\1_S$
and of the type
$\varepsilon (x_1')\ldots \varepsilon (x_k')\tau (y_1)\ldots
\tau (y_m)\1_{\Lambda}$, respectively, where $x_1', \ldots , x_k'\in V_+'$,
$y_1,\ldots y_m \in V_-$. Let $Deg \varepsilon (x') = Deg  \gamma (x') = 1$,
and
$Deg \tau (x) = Deg \beta (x) = -1$. Correspondingly, we obtain  $\Z$-gradings on the spaces of semi-infinite power modules: 
$S^{{\infty \over 2} + *}(V) = \oplus_{i\in\Z}S^{{\infty \over 2} + i}(V)$,
$\Lambda^{{\infty \over 2} + *}(V) = \oplus_{i\in\Z}
\Lambda^{{\infty \over 2} + i}(V)$.

Let $\lbrace e_i\rbrace _{i\in\Z}$ be a  {\it homogeneous} basis of $V$
 so that if
$i \in\Z$, then $e_i \in V_n$ for some $n\in \Z$, and if
$e_i\in V_n$, then $e_{i+1}\in V_n$ or $e_{i+1}\in V_{n+1}$.
Let $\lbrace e_i'\rbrace_{i\in\Z}$ be the dual basis.
Let $i_0 \in\Z$ be such that
$e_{i_0}\in V_K$ and $e_{i_0+1}\in V_{K+1}$.

Notice that one can think of  
$\Lambda ^{{\infty \over 2} + *}(V)$ as the vector
 space spanned by the elements
$w = e_{i_1}'\wedge e_{i_2}'\wedge\ldots$ such that there
 exists $N(w)\in\Z $
 such that $i_{n+1} = i_{n} - 1$ for $n > N(w)$.
Then $\1_{\Lambda} =  e_{i_0}'\wedge e_{i_0-1}'\wedge\ldots $
is a vacuum vector in this space. 
The actions of $\varepsilon (x'), \tau (x)$ are, respectively,
the exterior  multiplication and contraction in the space of 
 semi-infinite exterior products.

Let $\g = \oplus_{n \in\Z}\g_n$ be a graded Lie algebra over $\C$,
such that $dim \g_n < \infty$. 
Let $\phi$ be a representation of $\g$ in $V$ so that
$$\phi(\g_n)V_k \subset V_{k+n}.\eqno (2.1)$$
One can define the projective representations $\rho$ and $\pi$ of $\g$
in ${\Lambda}^{{\infty \over 2} + *}(V)$ and
$S^{{\infty \over 2} + *}(V)$, respectively:
$$\eqalignno{
&\rho(x) = \sum_{i\in\Z} : \tau(\phi(x)e_i)\varepsilon(e_i') :, &(2.2)\cr
&\pi(x) = \sum_{i\in\Z} :\beta (\phi(x)e_i)\gamma (e_i') :,&(2.3)\cr
}$$
where $x\in\g$,  and  where the double colons
$:\hbox { } :$ denote a normal ordering operation:
$$: \tau (e_j)\varepsilon (e_i') : =
\left\{
\eqalign{\tau (e_j)\varepsilon (e_i') \hbox { if } i \leq i_0\cr
- \varepsilon (e_i')\tau (e_j) \hbox { if } i > i_0\cr}\right\},\qquad
: \beta (e_j)\gamma (e_i') : =
\left\{
\eqalign{\beta (e_j)\gamma (e_i') \hbox { if } i \leq i_0\cr
 \gamma (e_i')\beta (e_j) \hbox { if } i > i_0\cr}\right\}.\eqno(2.4)
$$
Thus
$$\rho (x)\1_{\Lambda} = \pi (x)\1_S = 0 \hbox{ for } x \in \g_0\eqno
(2.5)$$
and
$$\eqalignno{
&[\rho (x), \rho (y)] = \rho ([x, y]) + c_{\Lambda}(x, y), &(2.6)\cr
&[\pi (x), \pi (y)] = \pi ([x, y]) + c_S(x, y), \cr
}$$
where $x, y \in \g$ and $c_{\Lambda}, c_S$ are 2-cocycles. Notice that
$c_{\Lambda} = - c_S$.
Let
$$W^{{\infty \over 2} + *}(V) = 
S^{{\infty \over 2} + *}(V)\otimes
\Lambda^{{\infty \over 2} + *}(V). \eqno(2.7)$$
Since the cocycles corresponding to the projective representations
cancel, the representation $\theta (x) = \rho (x) + \pi (x)$
of $\g$ in $W^{{\infty \over 2} + *}(V)$ is well-defined.
We define a $\Z$-grading on $W^{{\infty \over 2} + *}(V)$
setting
$$W^{{\infty \over 2} + i}(V) =  \oplus_{2l+j=i}
S^{{\infty \over 2} + l}(V)\otimes
\Lambda^{{\infty \over 2} + j}(V). \eqno(2.8)$$

Let $V = \g = \oplus_{n\in\Z}\g_n$ and $\phi$ be the adjoint 
representation of $\g$.
We define two differentials on the space $W^{{\infty \over 2} + *}(\g)$:
$$\eqalignno{
&d = \sum_{i<j}:\tau([e_i, e_j])\varepsilon (e_j')\varepsilon (e_i'): +
\sum_{i,j}:\beta([e_j, e_i])\gamma (e_i')\varepsilon (e_j'):, &(2.9)\cr
&\h1 = \sum_i \gamma (e_i')\tau (e_i).\cr
}$$
We obtain  {\it the semi-infinite Weil complex}
$$\lbrace W^{{\infty \over 2} + *}(\g), \hbox{ }d + \h1 \rbrace,
 \eqno
(2.10)$$
The differential $d$ is
the analogue of the classical differential for the Lie algebra
\hfil\break
(co)homology,
and $\h1$ is the analogue of the Koszul differential. Notice that
$$d^2 = 0, \h1^2 = 0, [d, \h1] = 0,  (d + \h1)^2 = 0.
\eqno (2.11)$$
Notice also that if $\g$ is a finite-dimensional Lie algebra, 
then applying the definitions given above to the 
polarization $\g = \g_+ \oplus \g_-$, where
$\g_+ = \g$, $\g_- = 0$, we obtain  the classical Weil complex.

As in the case of the classical Weil complex, one can construct
two filtrations, $F^*_1$ and $F^*_2$,
on $W^{{\infty \over 2} + *}(\g)$:
$$F^p_1 = \oplus_{l+j\geq p} S^{{\infty \over 2} + l}(\g)\otimes
\Lambda^{{\infty \over 2} + j}(\g),\quad
F^p_2 = \oplus_{2l\geq p} S^{{\infty \over 2} + l}(\g)\otimes
\Lambda^{{\infty \over 2} + *}(\g).\eqno (2.12)$$
For filtration $F^*_1$ the complex is acyclic, the second term of the 
spectral sequence associated to  filtration $F^*_2$ is 
{\it the semi-infinite cohomology} of Lie algebra $\g$ with coefficients in its ``adjoint semi-infinite symmetric powers"
$H^{{\infty \over 2} + *}(\g,    S^{{\infty \over 2} + *}(\g))$ (see [FF]). 
Let 
$$W^{{\infty \over 2} + *}_{rel}(V) = \lbrace w\in 
W^{{\infty \over 2} + *}(V)\mid \tau (x)w = 0 \hbox{ for all } x\in V_0,
\theta (x)w = 0 \hbox{ for all } x\in \g_0
\rbrace. \eqno (2.13)$$
The differential $d$ preserves the space
$W^{{\infty \over 2} + *}_{rel}(\g)$ since
$$[d, \tau(x)] = d\tau(x) + \tau(x)d = \theta (x), \eqno
(2.14)$$
and
$$[d, \theta(x)] = 0. \eqno(2.15)$$
for any $x \in \g$.
The complex $\lbrace W^{{\infty \over 2} + *}_{rel}(\g), d\rbrace$
is called {\it the relative semi-infinite Weil complex}.
Its cohomology is called {\it the relative semi-infinite cohomology}
$H^{{\infty \over 2} + *}(\g, \g_0,    S^{{\infty \over 2} + *}(\g))$.

We fix $K = 0$ from this point on. Correspondingly, $V = V_+\oplus V_-$, 
where $V_+ = \oplus_{n>0}V_n$,
$V_- = \oplus_{n\leq 0}V_n$. 
  
\vskip 0.2in
\centerline{\tenbf 3. The $N = 2$ superconformal algebra}
\vskip 0.2in

Recall that the   $N = 2$ SCA is spanned by the
Virasoro
generators $\L_n$, Heisenberg generators $H_n$, two fermionic fields
$G^{\pm}_r$, and a central element $\centre$, where
$n  \in \Z, r \in  \Z + 1/2$, and where the non-vanishing commutation
relations are as follows, [FST]:
$$\eqalignno{
&[\L_n, \L_m] = (n - m) \L_{n+m} + {\centre\over 12}(n^3 - n)\delta _{n,
-m}, &(3.1)
\cr
&[\L_n, H_m] = -m H_{n+m}, 
[\L_n, G^{\pm}_r] = ({n\over 2} - r)G^{\pm}_{n+r}, \cr
&[G^{+}_r, G^{-}_s] = 2\L_{r+s} + (r - s)H_{r+s} + {\centre \over 3}
(r^2 - {1\over 4})
\delta _{r, -s}, \cr
&[H_n, H_m] = {\centre \over 3}n \delta _{n, -m}, 
[H_n, G^{\pm}_r] = \pm G^{\pm}_{n+r}. \cr
}$$
Let $Witt = \oplus_{i\in \Z}\C L_i$ be the Witt algebra:
$$[L_i, L_j] = (i - j)L_{i+j}. \eqno (3.2)$$
Let $\lambda, \mu \in \C$. Let
$\F1_{\lambda, \mu} = \oplus_{m\in \Z}\C u_m$ be a module over
$Witt$ defined as follows:
$$\phi (L_n)u_m = (-m + \mu -(n - 1)\lambda)u_{n+m}.\eqno (3.3)$$

{\it Remark 3.1.}-
The module $\F1_{\lambda, \mu} = \oplus_{m\in \Z}\C u_m$ 
is isomorphic to the module 
$\F1_{-\lambda, \mu + 1} = \oplus_{j\in\Z}\C f_j$
over the Witt algebra
defined in [Fu]. The isomorphism is given
by the correspondence $u_m \leftrightarrow f_{-m-1}$.

{\sl Theorem 3.1.}-  
The space $W^{{\infty \over 2} + *}(\F1_{\lambda, \mu})$
is a module over the $N = 2$
SCA with central charge $3 - 6\lambda$.

{\it Proof.}- Set
$$\h1_n = {1\over {\sqrt 2}}G^{+}_{n-{1\over 2}},
\p1_n = {1\over {\sqrt 2}}G^{-}_{n+{1\over 2}}. \eqno (3.4)$$
 We define a representation  of $Witt$ in 
  $W^{{\infty \over 2} + *}(\F1_{\lambda, \mu})$ as follows:
$$\theta (L_n) = \sum_{m\in\Z} (-m + \mu -n\lambda + \lambda) \bigl(:\tau(u_{m+n})
\varepsilon(u_{m}'): + :\beta(u_{m+n})\gamma (u_{m}'):).\eqno (3.5)$$
Let us extend $\theta$ to a representation  of the $N = 2$ SCA in
$W^{{\infty \over 2} + *}(\F1_{\lambda, \mu})$:
$$\eqalignno{
&\theta(H_n) = \lambda  \sum_{m\in\Z} :\tau(u_m)\varepsilon(u_{m+n}'): +
(\lambda - 1) \sum_{m\in\Z} : \beta(u_m)\gamma(u_{m+n}'): + \mu\delta_{n,0},
&(3.6)\cr
&\theta(\h1_n) = \sum_{m\in\Z} \gamma(u_{m+n}')\tau(u_m),\cr
&\theta(\p1_n) = \sum_{m\in\Z}(m - \mu - (n +
1)\lambda)\beta(u_{m-n})\varepsilon(u_{m}'),
\cr
&\theta(\L_n) = -\theta(L_{-n}) + {{n+1}\over 2}\theta(H_n).\cr
}$$
We calculate the central charge by
checking the commutation relations on the vacuum vector
$\1 = \1_S\otimes \1_{\Lambda}$.
Let $n > 0$. Then
$$\eqalignno{
&\theta([H_n, H_{-n}])\1 = -\theta(H_{-n})\theta(H_n)\1 = 
 - \theta(H_{-n})\bigl(\lambda \sum_{m=1-n}^0\tau(u_m)\varepsilon(u_{m+n}')
+ &(3.7)\cr
&(\lambda - 1)\sum_{m=1-n}^0\beta(u_m)\gamma(u_{m+n}')\bigr)\1 = 
-\lambda^2\sum_{m=1-n}^0\tau(u_{m+n})\varepsilon(u_{m}')
\tau(u_m)\varepsilon(u_{m+n}')\1 - \cr
&(\lambda - 1)^2
\sum_{m=1-n}^0\beta(u_{m+n})\gamma(u_{m}')
\beta(u_m)\gamma(u_{m+n}')\1 = 
\bigl(-\lambda^2n - (\lambda - 1)^2(-n)\bigr)\1 = n(1 - 2\lambda)\1,\cr
}$$
since $\varepsilon(u_i')\tau(u_i) + \tau(u_i)\varepsilon(u_i') = 1$, and
$\gamma(u_i')\beta(u_i) - \beta(u_i)\gamma(u_i') = 1$.
Hence,
$$\theta([H_n, H_m])\1 = n(1 - 2\lambda)\delta_{n, -m}\1. \eqno (3.8)$$
Thus the central charge is    $3 - 6\lambda$.
The other commutation relations on the vacuum vector $\1$ are calculated in
the same way.
$$\eqno \SQ$$

{\it Remark 3.2.}- In  the case when
$\lambda = - 1, \mu = 1$, the module $\F1_{\lambda, \mu}$ is the
adjoint representation of $Witt$. Thus we obtain a representation
of the $N = 2$ SCA in the semi-infinite Weil complex
of the Witt algebra (cf. [G]).

{\sl Theorem 3.2.} 
Let $V$ be a complex finite-dimensional vector space,
$\tilde{V} = V\otimes \C[t, t^{-1}]$. 
There exists a representation of the $N = 2$  SCA
in $W^{{\infty \over 2} + *}(\tilde{V})$  
with central charge $3dim V$.

{\it Proof.}- There is the natural $\Z$-grading
$\tilde{V} = \oplus_{n\in\Z}\tilde{V}_n$, where 
$\tilde{V}_n = V\otimes  t^n$.
Let  $u$ run through a fixed basis of $V$,
$u_n$ stand for $u\otimes t^n$, and let
$\lbrace u_n'\rbrace$ be the dual basis of $\tilde{V}'$.
Define the following quadratic expansions by analogy with (3.5) and (3.6),
where $\lambda = 0, \mu = 0$:
$$\eqalignno{
&L_n = -\sum_u\sum_{m\in\Z} \bigl( m : \tau (u_{m+n})\varepsilon (u_m') :
+ m : \beta (u_{m+n})\gamma (u_m') :\bigr)\cr
&H_n = - \sum_u\sum_{m\in\Z} : \gamma (u_{m+n}')\beta (u_m) :\cr
&\h1_n = \sum_u\sum_{m\in\Z} \gamma (u_{m+n}')\tau (u_m), &(3.9)\cr
&\p1_n = \sum_u\sum_{m\in\Z} m\beta (u_{m-n})\varepsilon (u_m').\cr
}$$
Set
$$\L_n = -L_{-n} + {{n+1} \over 2}H_n.\eqno (3.10)$$
Then $\L_n, H_n, \h1_n$, and $\p1_n$ span the centerless
$N = 2$  SCA. 

Let $n > 0$. Then $H_{-n}\1 = 0$. Hence
$$\eqalignno{
&[H_n, H_{-n}]\1 = -H_{-n}\bigl(-\sum_u\sum_{m=1-n}^{0}
\gamma(u_{m+n}')\beta(u_m)\bigr)\1 = \cr
&\bigl(-\sum_u\sum_{m=1}^{n}\gamma(u_{m-n}')\beta(u_m)\bigr)
\bigl(\sum_u\sum_{m=1-n}^{0}\gamma(u_{m+n}')\beta(u_m)\bigr)\1= &(3.11)\cr
&-\sum_u\sum_{m=1-n}^{0}\gamma(u_{m}')\beta(u_{m+n})
\gamma(u_{m+n}')\beta(u_{m})\1 = -dim V (-n)\1,\cr
}$$  
since $\gamma(u_{i}')\beta(u_i) - \beta(u_i)\gamma(u_{i}') = 1$.
Notice that 
$$[H_n, H_{m}]\1 = 0, \hbox { if } m\not= -n.\eqno (3.12)$$
Hence
$$[H_n, H_m]\1 = ndim V \delta_{n, -m}\1. \eqno(3.13)$$
Thus the central charge is $3dim V$.
$$\eqno \SQ$$

{\sl Corollary 3.1.}- 
Let $\g$ be a complex finite-dimensional Lie algebra, let
$\tilde{\g} = \g\otimes \C[t, t^{-1}]$.
There exists a representation of the $N = 2$  SCA in
$H^{{\infty \over 2} + *}(\tilde{\g},    S^{{\infty \over 2} + *}(\tilde{\g}))$
with central charge $3dim \g$.

{\it Proof.}- We will show that the  expansions 
(3.9) commute with the differential $d$.
Recall that
$$d = d^{(1)} + d^{(2)},\eqno (3.14)$$
where
$$\eqalignno{
&d^{(1)} = (1/2)\sum_{u,v,i,j} : \tau ([u_i, v_j])\varepsilon (v_j')
\varepsilon (u_i') :, &(3.15)\cr
&d^{(2)} = \sum_{u,v,i,j} 
: \beta ([u_i, v_j])\gamma (v_j')  \varepsilon (u_i') :, \cr
}$$
$u, v$ run through a fixed basis of $\g$, and $i, j\in\Z$.
Then
$$\eqalignno{
&[L_n, d^{(1)}] = (1/2)\sum_{u,v,i,j} : -(i + j) \tau ([u, v]_{i+j+n})
\varepsilon (v_j')\varepsilon (u_i') : + \cr
&: \tau ([u_i, v_j])(j - n)\varepsilon (v_{j-n}')\varepsilon (u_i') :
+ : \tau ([u_i, v_j])\varepsilon (v_{j}')(i - n)\varepsilon (u_{i-n}') :
  = 0 &(3.16)\cr
}$$
and
$$\eqalignno{
&[L_n, d^{(2)}] = 
\sum_{u,v,i,j} : -(i + j)\beta ([u, v]_{i+j+n})\gamma (v_j')  \varepsilon (u_i') : + \cr
&: \beta ([u_i, v_j])(j - n)\gamma (v_{j-n}')
\varepsilon (u_i') : + : \beta ([u_i, v_j])\gamma (v_j')(i - n)
\varepsilon (u_{i-n}') : = 0. &(3.17)\cr
}$$
Clearly,
$$[H_n, d^{(1)}] = 0,\eqno (3.18)$$
and
$$[H_n, d^{(2)}] = \sum_{u,v,i,j} : -\beta ([u, v]_{i+j-n})\gamma (v_j')
\varepsilon (u_i') : + 
\beta ([u, v]_{i+j})\gamma (v_{j+n}')\varepsilon (u_i') : = 0.\eqno (3.19)$$
Next,
$$\eqalignno{
&[\h1_n, d^{(1)}] = (1/2)\sum_{u,v,i,j} : - \tau ([u, v]_{i+j})
\gamma (v_{j+n}')\varepsilon (u_i') : + 
: \tau ([u, v]_{i+j})
\varepsilon (v_{j}')\gamma (u_{i=n}') : = \cr
&-\sum_{u,v,i,j} : \tau ([u, v]_{i+j})\gamma (v_{j+n}')\varepsilon (u_i') :,
&(3.20)\cr
&[\h1_n, d^{(2)}] = \sum_{u,v,i,j} : \tau ([u, v]_{i+j-n})
\gamma (v_{j}')\varepsilon (u_i') : + 
: \beta ([u, v]_{i+j})
\gamma (v_{j}')\gamma (u_{i+n}') : = \cr
&\sum_{u,v,i,j} : \tau ([u, v]_{i+j})
\gamma (v_{j+n}')\varepsilon (u_i') :, &(3.21)\cr
}$$
since $\sum_{u,v,i,j} : \beta ([u, v]_{i+j})
\gamma (v_{j}')\gamma (u_{i+n}') : = 0$. Hence 
$$[\h1_n, d^{(2)}] = - [\h1_n, d^{(1)}]. \eqno (3.22)$$
Finally,
$$\eqalignno{
&[\p1_n, d^{(1)}] = (1/2)\sum_{u,v,i,j} : (i + j)\beta ([u, v]_{i+j-n})
\varepsilon (v_j')\varepsilon (u_i') :, &(3.23)\cr
&[\p1_n, d^{(2)}] = \sum_{u,v,i,j} : -\beta ([u, v]_{i+j})(j + n)
\varepsilon (v_{j+n}')\varepsilon (u_i') : = &(3.24)\cr
&\sum_{u,v,i,j} : -\beta ([u, v]_{i+j-n})j
\varepsilon (v_{j}')\varepsilon (u_i') : = 
-(1/2)\sum_{u,v,i,j} : (j + i)\beta ([u, v]_{i+j-n})
\varepsilon (v_{j}')\varepsilon (u_i') : .\cr
}$$ 
Hence 
$$[\p1_n, d^{(2)}] = - [\p1_n, d^{(1)}]. \eqno (3.25)$$
$$\eqno \SQ$$

\vskip 0.2in
\centerline{\tenbf 4. The superconformal algebras $S'(2, \alpha)$}
\vskip 0.2in

Recall the necessary definitions, [KL]. Let
$W(N)$ be the superalgebra of all derivations of
$\C [t, t^{-1}]\otimes \Lambda (N),$
where $\Lambda(N)$ is the Grassmann algebra in $N$ variables
$\theta_1, \ldots, \theta_N$, and $p(t) = \bar{0}$, $p(\theta_i) =
\bar{1}$ for $i = 1, \ldots, N$.
Let $\partial_i$ stand for $\partial/{\partial\theta_i}$, and
$\partial_t$ stand for $\partial/{\partial t}$.
Let
$$S(N, \alpha) = \lbrace D \in W(N) \mid
Div(t^{\alpha}D) = 0\rbrace \hbox { for } \alpha \in \C.\eqno (4.1)$$
Recall that
$$Div(f\partial_t + \sum_{i=1}^N f_i \partial_i) =
\partial_tf + \sum_{i=1}^N (-1)^{p(f_i)}
\partial_i f_i\eqno (4.2)$$
where $f, f_i \in \C [t, t^{-1}]\otimes \Lambda (N)$, and
$$Div(fD) = Df + fDiv D, \eqno (4.3)$$
where $f$ is an even function.
Let $S'(N, \alpha) =  [S(N, \alpha), S(N, \alpha)]$ be the derived
superalgebra.
Assume that $N > 1$. If ${\alpha}\not\in\Z$, then
$S(N, \alpha)$ is simple, and if 
${\alpha} \in \Z$, then
$S'(N, \alpha) $ is a  simple ideal of $S(N, \alpha)$ of codimension 1:
$$0\rightarrow S'(N, \alpha)\rightarrow
 S(N, \alpha)\rightarrow \C t^{-\alpha} \theta_1 \cdots  \theta_N
\partial_t\rightarrow 0.\eqno (4.4)$$
Notice that
$$S(N, \alpha)\cong S(N, \alpha + n) \hbox{ for } n\in\Z. \eqno (4.5)$$
The superalgebra $S'(N, \alpha)$ has, up to equivalence, only one
non-trivial cocycle if and only if
$N = 2$, which is important for our task.
Let $$\lbrace \L_n^{\alpha}, E_n, H_n,
 F_n, \h1_n^{\alpha}, \p1_n,
 \x_n, \y_n^{\alpha}\rbrace _{n \in\Z}
\eqno (4.6)$$
be  the basis of $S'(2,\alpha)$ defined as follows:
$$\eqalignno{
&\L_n^{\alpha} = -t^n(t \partial_t + {1\over 2}(n + \alpha + 1)
(\theta_1 \partial_1 +
\theta_2 \partial_2 )), &(4.7)\cr
& E_n = t^n \theta_2 \partial_1, \cr
& H_n = t^n (\theta_2 \partial_2 -
\theta_1 \partial_1), \cr
&F_n = t^n \theta_1 \partial_2, \cr
&\h1_n^{\alpha} = t^{n}\theta_2\partial_t -
(n + \alpha)t^{n-1}\theta_1\theta_2\partial_1,\cr
&\p1_n = -t^{n+1}\partial_2,\cr
&\x_n = t^{n+1}\partial_1, \cr
&\y_n^{\alpha} = t^{n}\theta_1\partial_t +
(n + \alpha)t^{n-1}\theta_1\theta_2\partial_2.\cr
}$$
The non-vanishing commutation relations between these elements
 are
$$\eqalignno{
&[\L_n^{\alpha}, \L_k^{\alpha}] = (n - k)\L_{n+k}^{\alpha}, &(4.8)\cr
&[E_n, F_k] = H_{n+k}, [H_n, E_k] = 2E_{n+k}, [H_n, F_k] = -2F_{n+k},\cr
&[\L_n^{\alpha} , E_k] = -kE_{n+k},
 [\L_n^{\alpha}, H_k] = -kH_{n+k},
[\L_n^{\alpha}, F_k] = -kF_{n+k},\cr
&[\L_n^{\alpha}, \h1_k^{\alpha}] =
 {1\over 2}(n-2k+1-\alpha)\h1_{n+k}^{\alpha},
[\L_n^{\alpha}, \p1_k] = {1\over 2}(n-2k-1+\alpha)\p1_{n+k},\cr
&[\L_n^{\alpha}, \x_k] =
{1\over 2}(n-2k-1+\alpha)\x_{n+k},
[\L_n^{\alpha}, \y_k^{\alpha}] = {1\over 2}(n-2k+1-\alpha)
\y_{n+k}^{\alpha}, \cr
&[E_n, \y_k^{\alpha}] = \h1_{n+k}^{\alpha},
 [F_n, \h1_k^{\alpha}] = \y_{n+k}^{\alpha},
[E_n, \p1_k] = \x_{n+k},
 [F_n, \x_k] = \p1_{n+k},\cr
&[H_n, \h1_k^{\alpha}] = \h1_{n+k}^{\alpha},
 [H_n, \y_k^{\alpha}] = -\y_{n+k}^{\alpha},
[H_n, \x_k] = \x_{n+k},
 [H_n, \p1_k] = -\p1_{n+k}, \cr
&[\h1_n^{\alpha}, \x_k] =
 (k+1-n-\alpha)E_{n+k},
[\p1_n, \y_k^{\alpha}] = (k-n-1+\alpha)F_{n+k},\cr
&[\h1_n^{\alpha}, \p1_k] = \L_{n+k}^{\alpha}
 - {1\over 2}(k - n + 1 - \alpha)H_{n+k},\cr
 &[\x_n, \y_k^{\alpha}] = -\L_{n+k}^{\alpha} +
 {1\over 2}(k - n - 1 + \alpha)H_{n+k}.\cr
}$$
A non-trivial cocycle on ${S}'(2, \alpha)$ is
$$\eqalignno{
&c(\L_n^{\alpha}, \L_k^{\alpha}) = {\centre\over 12}n(n^2 - 1)\delta_{n, -k},
&(4.9)\cr
&c(E_n, F_k) = {\centre\over 6}n\delta_{n,-k},
c(H_n, H_k) = {\centre \over 3}n\delta_{n, -k}, \cr
&c(\h1_n^{\alpha}, \p1_k) = 
{\centre\over 6}\bigl((n - 1 + {{\alpha + 1}\over 2})^2 -
 {1\over 4}\bigr)\delta_{n,-k},\cr
&c(\x_n, \y_k^{\alpha}) =
- {\centre\over 6}
\bigl((-n - 1 + {{\alpha + 1}\over 2})^2 - {1\over 4}\bigr)\delta_{n,-k};\cr
}$$
see [KL]. Let $\hat{S}'(2, \alpha)$ be the corresponding central
extension of $S'(2, \alpha)$. 
In particular,  $\hat{S}'(2, 0)$ is isomorphic to the 
$N = 4$  SCA (see [Ad]).

{\it Remark 4.1.}- Notice that
$$\eqalignno{
&S'(2, \alpha)_{\bar 0} = Witt \+ \tilde{\s\l}(2), \hbox{ where }  &(4.10)\cr
&Witt = \langle \L_{n}^{\alpha} \rangle_{n\in\Z},
\tilde{\s\l}(2) = \langle E_{n}, H_{n}, F_{n}\rangle_{n\in\Z},\cr
}$$
and
$$S'(2, \alpha)_{\bar 1} = 
\langle \h1_n^{\alpha},\y_n^{\alpha} \rangle_{n\in\Z}\oplus
\langle \p1_n, \x_n \rangle_{n\in\Z} \eqno (4.11)$$
is a direct sum of two standard (odd) $\tilde{\s\l}(2)$-modules. 
 
{\it Remark 4.2.}- For any $\alpha \in \C$ one can consider
the subalgebra of $\hat{S}'(2, \alpha)$,
spanned by $\L^{\alpha}_n, H_n,\h1_n^{\alpha},
\p1_n$, and $\centre$.
Thus we obtain a one-parameter family of superalgebras, 
which are isomorphic to the $N = 2$  SCA.
The isomorphism
$$\varphi : \langle \L^{\alpha}_n, H_n,\h1_n^{\alpha},
\p1_n, \centre \rangle \longrightarrow  \langle \L_n, H_n,\h1_n,
\p1_n, \centre \rangle\eqno (4.12)$$
is given as follows:
$$\eqalignno{
&\varphi (\L_n^{\alpha}) =  \L_n - {\alpha\over 2}H_n +
{\alpha^2\over 24}\delta_{n,0}\centre, &(4.13)\cr
&\varphi (H_n) = H_n - {\alpha\over 6}\delta_{n,0}\centre,\cr
&\varphi (\h1_n^{\alpha}) = \h1_n,\
\varphi (\p1_n) = \p1_n, \varphi (\centre) = \centre.\cr
}$$
Notice that formulae (4.13) correspond
to the spectral flow transformation for the $N = 2$  SCA
(cf. [FST]).

Let $Der S'(2,\alpha)$ be the Lie superalgebra
 of all derivations of $S'(2,\alpha)$,
and 
\hfil\break
$Der_{ext} S'(2,\alpha)$ be the exterior
derivations of $S'(2,\alpha)$ (see [Fu]).

{\sl Theorem 4.1.}- 
\hfil\break
1) If $\alpha \in \Z$, then 
$Der_{ext} S'(2,\alpha) \cong 
\S\L(2) = \langle \E1,\H1,\F1 \rangle$, where
$$[\H1, \E1] = 2\E1, 
[\H1, \F1] = - 2\F1,
[\E1, \F1] = \H1. \eqno (4.14)$$
The action of $\S\L(2)$  is given as follows:
$$\eqalignno{
&[\E1, \h1_k^{\alpha}] = \x_{k-1+\alpha},
[\E1, \y_k^{\alpha}] = \p1_{k-1+\alpha}; &(4.15)\cr
&[\F1, \x_k] = \h1_{k+1-\alpha}^{\alpha},
[\F1, \p1_k] =  \y_{k+1-\alpha}^{\alpha}; \cr
&[\H1, \x_k] = \x_k, [\H1, \h1_k^{\alpha}] =
 - \h1_k^{\alpha}, \cr
&[\H1, \p1_k] = \p1_k, [\H1, \y_k^{\alpha}] =
 - \y_k^{\alpha}.  \cr
}$$
2) If $\alpha \in \C\setminus \Z$, then
$Der_{ext}S'(2,\alpha) = \langle \H1 \rangle$.

{\it Proof.}-
Recall that the exterior derivations of a Lie (super) algebra
can be identified with its first cohomology with coefficients
in the adjoint representation (see [Fu]). Thus 
$$Der_{ext} S'(2,\alpha) \cong H^1(S'(2,\alpha), S'(2,\alpha)).
 \eqno(4.16)$$
The superalgebra $S'(2,\alpha)$ has the following $\Z\pm \alpha$-grading $deg$:
$$\eqalignno{
&deg \L_n^{\alpha} =  n, deg E_n =  n + 1 - \alpha,
 deg F_n =  n - 1 + \alpha, deg H_n =  n,
&(4.17)\cr
&deg \h1_n^{\alpha} =  n, deg \p1_n =  n,
 deg \x_n =  n + 1 - \alpha, deg \y_n^{\alpha} =  n - 1 + \alpha.
\cr
}$$
Let
$$L_0 = - \L_{0}^{\alpha} + {1\over 2}(1 -  \alpha)H_{0}.
 \eqno (4.18)$$
 Then
 $$[L_0, s] = (deg s)s \eqno (4.19)$$
 for a homogeneous $s \in S'(2,\alpha)$.
Accordingly,
$$[L_0, D] = (deg D)D \eqno (4.20)$$
for a homogeneous   $D \in Der_{ext} S'(2,\alpha)$.
On the other hand, since the action of a Lie superalgebra
on its cohomology is trivial (see [Fu]), then one must have
$$[L_0, D] = 0.\eqno (4.21)$$
Hence the non-zero elements of $Der_{ext} S'(2,\alpha)$ have $deg = 0$,
and they preserve the superalgebra
$S'(2,\alpha)_{deg=0}$. Let $\alpha\in\Z$. Then
one can check that the exterior derivations of  $S'(2,\alpha)_{deg=0}$
form an $\s\l(2)$, and  extend
them to the exterior derivations of $S'(2,\alpha)$ as in (4.15). 
One should also note that
if the restriction of a derivation of $S'(2,\alpha)$ to $S'(2,\alpha)_{deg=0}$
is zero, then this derivation is inner.

Finally, notice that the exterior derivations
$\E1$ and $\F1$ interchange $\lbrace\h1_k^{\alpha}\rbrace$ with
$\lbrace\x_k\rbrace$.
If $\alpha \not\in\Z$, then
$deg \h1_k^{\alpha} - deg \x_n \not\in \Z$  for any
$k,n \in\Z$. Hence  $\E1$ and $\F1$ cannot have $deg = 0$.
By this reason, $Der_{ext} S'(2,\alpha) = \langle \H1\rangle$
for $\alpha \in\C\setminus \Z$.
$$\eqno \SQ$$

{\it Remark 4.3.}- If $\alpha\in\Z$, then
one can identify $\F1$ with 
$-t^{-\alpha} \theta_1\theta_2\partial_t$ (see (4.4)).
\vskip 0.2 in
\centerline {\tenbf 5. An action  of
$\hat{S}'(2, \alpha)$ on the  semi-infinite Weil
complex of a loop algebra}
\vskip 0.2 in
\noindent
We will consider a more general case, i. e. when
$V$ is a complex
 finite-dimensional vector space, and
$\tilde{V} = V\otimes \C[t, t^{-1}]$. 
Let $\hat{D}er S'(2,\alpha)$ be a non-trivial central extension of
${D}er S'(2,\alpha)$.

{\sl Theorem 5.1.}-
\hfil\break
1) The space
$W^{{\infty \over 2} + *}(\tilde{V})$, where $\alpha \in \C$, is a module over
 $\hat{S}'(2,\alpha)$ with central charge $3dim V$;
\hfil\break
2) if $\alpha\in\C\setminus\Z$, then $W^{{\infty \over 2} + *}(\tilde{V})$ is a
module over  $\hat{D}er S'(2,\alpha)$.
\hfil\break
{\it Proof.} 
Let  $u$ run through a fixed basis of $V$,
$u_n$ stand for $u\otimes t^n$, and
$\lbrace u_n'\rbrace$ be the dual basis of $\tilde{V}'$.
One can define 
a representation of $Witt$
in $W^{{\infty \over 2} +
*}(\tilde{V})$ by analogy with (3.5), where $\lambda = 0, \mu = \alpha/2$:
$$\theta (L_n) =
- \sum_u\sum_{m}(m - {\alpha\over 2})
\bigl(: \tau (u_{m+n})\varepsilon (u_{m}') : +
 : \beta (u_{m+n})\gamma (u_{m}') :\bigr), \eqno (5.1)$$
then extend it to a representation of
the $N = 2$  SCA, and apply 
 (4.13). We obtain 
the following representation of $\hat{S}'(2, \alpha)$:
$$\eqalignno{
&\theta (H_n) =
- \sum_u\sum_{m} :  \beta (u_m)\gamma (u_{m+n}'), &(5.2)\cr
&\theta (\L_n^{\alpha}) =
-\theta (L_{-n}) +
 {{n + 1 - \alpha}\over 2}\theta (H_n) + ({\alpha \over 4} -
{\alpha^2\over 8})dimV\delta_{n,0}, \cr
 &\theta(\h1_n^{\alpha}) =
\sum_u\sum_{m} \gamma (u_{m+n}')\tau (u_m), \cr
&\theta(\p1_n) = \sum_u\sum_{m} (m - {\alpha\over 2})
\beta (u_{m-n})\varepsilon (u_m'), \cr
&\theta(E_n) =
 -(1/2)i\sum_u\sum_{m}\gamma(u_m')\gamma(u_{1-m+n}'), \cr
&\theta(F_n) =
 -(1/2)i\sum_u\sum_{m}\beta(u_m)\beta(u_{1-m-n}), \cr
&\theta(\y_n^{\alpha}) = i\sum_u\sum_{m}
 \beta (u_{m})\tau (u_{1-m-n}), \cr
&\theta(\x_n) = -i\sum_u\sum_{m}(m - {\alpha\over 2})
 \gamma (u_{1-m+n}')\varepsilon (u_m'), \cr
&\theta (\H1) = - \sum_u\sum_{m} : \tau(u_m)\varepsilon(u_{m}'). \cr 
}$$
One can check that the central charge is $3dim V$ in the same way as in
Theorem 3.2.
$$\eqno \SQ$$

{\sl Theorem 5.2.}-
Let $\g$ be a complex finite-dimensional Lie algebra
endowed with a non-degenerate invariant symmetric
bilinear form. 
Then 
$H^{{\infty \over 2} + *}(\tilde{\g},    S^{{\infty \over 2} + *}(\tilde{\g}))$
is a module over $\hat{S}'(2, 0)$
with central charge $3dim \g$.

{\it  Proof.}-
Let $\lbrace v_i\rbrace$ be a basis of $\g$
so that with respect to the given form
$\langle v_i, v_j\rangle = \delta_{i,j}$. Let $u$ run through this basis.
Then by Theorem  5.1, there is a representation of
$\hat{S}'(2, 0)$ in $W^{{\infty \over 2} + *}(\tilde{\g})$.
Notice that we can  identify the elements
of $S'(2, 0)$ with the quadratic expansions obtained by putting 
$\alpha = 0$ in the equations (5.2).
One can check that the commutation relations (4.8) (where $\alpha = 0$)
are fulfilled.
One can notice that 
$$[S'(2,0), d] = 0. \eqno (5.3)$$
In fact, since $\langle\cdot, \cdot\rangle$ is an invariant symmetric bilinear form on
$\g$, then
the elements $E_n, H_n$, and $F_n$
commute with $\pi (g)$ for any $g\in\tilde{\g}$. Hence they commute with $d$.
According to Corollary 3.1,
$$[\h1_n^0, d] = [\p1_n, d] = 0.\eqno (5.4)$$
Recall that 
$$S'(2,0)_{\bar 1} = \langle \h1_n^0, \y_n^0, \p1_n, \x_n \rangle_{n\in\Z}. 
\eqno (5.5)$$
Since 
$$[E_n, \p1_k] = \x_{n+k}, [F_n, \h1_k^0] = \y_{n+k}^0,\eqno (5.6)$$
then
$$[S'(2,0)_{\bar 1}, d] = 0. \eqno (5.7)$$
Since
$$S'(2,0)_{\bar 0} = [S'(2,0)_{\bar 1}, S'(2,0)_{\bar 1}],
\eqno (5.8) $$
then (5.3) follows.
$$\eqno \SQ$$
To define an action of $\hat{D}er S'(2,0)$, 
one should consider a {\it relative} semi-infinite Weil complex.

Let $\g$ be a complex finite-dimensional Lie algebra, $\phi$ be a
representation of $\g$ in $V$,
$\langle\cdot,\cdot\rangle$ be a non-degenerate $\g$-invariant symmetric
bilinear form on $V$.
One can naturally extend $\phi$ to a   
representation of $\tilde{\g}$
in $\tilde{V}$:
$$\phi(g\otimes t^n)(v\otimes t^k) = (\phi(g)v) \otimes t^{n+k},
\hbox{ for } g\in\g, v\in V. \eqno (5.9)$$

{\sl Theorem 5.3.}- The space   
$W^{{\infty \over 2} + *}_{rel}(\tilde{V})$
is a module over
$\hat{D}er S'(2,0)$ with central charge $3dim V$.

{\it Proof.}-
Let $\lbrace v_i\rbrace$ be a basis of $V$
so that
$\langle v_i, v_j\rangle = \delta_{i,j}$. Let $u$ run through this basis.
Then by Theorem  5.1, there is a representation of
$\hat{S}'(2, 0)$ in $W^{{\infty \over 2} + *}(\tilde{V})$.
We can  identify the elements
of $S'(2, 0)$ with the  expansions (5.2) where $\alpha = 0$.

Since
the form $\langle\cdot ,\cdot \rangle$ is $\g$-invariant, then there is an action of  
$\hat{S}'(2, 0)$ on $W^{{\infty \over 2} + *}_{rel}(\tilde{V})$.
To extend this representation to $\hat{D}er S'(2,0)$, we have to define it on
$\S\L(2) =
\langle \F1, \H1, \E1 \rangle$.
Let 
$$\eqalignno{
&\E1 = i\sum_u\sum_{m>0}m\varepsilon(u_{-m}')\varepsilon(u_{m}'),&(5.10)\cr
&\H1 = - \sum_u\sum_{m\not=0} : \tau(u_m)\varepsilon(u_{m}') :, \cr
&\F1 = -i\sum_u\sum_{m>0}(1/m)\tau(u_{m})\tau(u_{-m}).\cr  
}$$
Notice that $\S\L(2)$ acts on
$W^{{\infty \over 2} + *}_{rel}(\tilde{V})$.
The commutation relations between $\E1, \H1, \F1$ and the elements
of  $S'(2,0)$
 coincide with the relations
(4.15), where $\alpha = 0$,  up to some terms which contain 
elements $\tau(u_0)$.
Since the action of  $\tau(u_0)$  on
$W^{{\infty \over 2} + *}_{rel}(\tilde{V})$ is trivial, then
a representation of $\hat{D}er S'(2,0)$ in $W^{{\infty \over 2} +
*}_{rel}(\tilde{V})$ is well-defined.
$$\eqno\SQ$$
 
{\sl Corollary 5.1.}-
$H^{{\infty \over 2} + *}
(\tilde{\g}, \tilde{\g}_0,   S^{{\infty \over 2} + *}(\tilde{\g}))$
is a module over $\hat{S}'(2, 0)$
with central charge $3dim \g$.

{\it Proof.}- Follows from Theorem 5.2.
$$\eqno\SQ$$
 \vskip 0.2in
\centerline{\tenbf 6. Relative semi-infinite cohomology and K\"ahler geometry}
\vskip 0.2in

Let $M$ be a compact  K\"ahler manifold with associated $(1, 1)$-form 
$\omega$, let 
$dim_{\C} M = n$.
 There exists a number of  operators on the space $A^*(M)$ of differential forms on $M$ such as $\partial, \bar{\partial}, d, d_c$, their corresponding adjoint operators and the associated
Laplacians (see [GH]). Recall that
$$\eqalignno{
&\partial: A^{p,q}(M) \rightarrow A^{p+1,q}(M), &(6.1)\cr
&\bar{\partial}: A^{p,q}(M) \rightarrow A^{p,q+1}(M),\cr
&d = \partial + \bar{\partial}, \cr
&d_c = i(\partial - \bar{\partial}), \cr
&\triangle = dd^* + d^*d = 2\triangle_{\partial} = 2\triangle_{\bar{\partial}}.\cr
}$$ 
The Hodge $\star$-operator maps
$$ \star : A^{p,q}(M)\longrightarrow A^{n-q,n-p}(M),\eqno (6.2)$$
so that $\star^2 = (-1)^{p+q}$ on  $A^{p,q}(M)$.
Correspondingly, the Hodge inner product is defined on  each of $A^{p,q}(M)$:
$$(\varphi, \psi) = \int_{M}\varphi\wedge \star\bar{\psi}.\eqno (6.3)$$
In addition, 
$A^*(M)$ admits an $\s\l(2)$-module structure.
Namely, 
$\s\l(2) = \langle L, H, \Lambda \rangle$, where
$$[L, \Lambda] = H, [H, L] = 2L, [H, \Lambda] = -2\Lambda. \eqno (6.4)$$
The operator
$$L : A^{p,q}(M)\rightarrow A^{p+1,q+1}(M), \eqno (6.5)$$
is defined by
$$L (\varphi) = \varphi \wedge \omega.\eqno (6.6)$$ 
Let $\Lambda = L^*$ be its adjoint operator:
$$\Lambda: A^{p,q}(M)\rightarrow A^{p-1,q-1}(M),\eqno (6.7)$$ 
and
$$H\mid_{A^{p,q}(M)} = p + q - n.\eqno (6.8)$$
According to the Lefschetz theorem, there exists the corresponding
action of $\s\l(2)$ on $H^*(M)$.
These operators satisfy a series of identities,
known as the Hodge identities (see  [GH]).
Consider the Lie superalgebra spanned by the classical operators:   
$$\S := \langle \triangle, L , H, \Lambda, 
d, d^*, d_c, d_c^*\rangle\eqno (6.9)$$
The non-vanishing commutation relations in $\S$ are as follows:
$$\eqalignno{
&[L, \Lambda] = H, [H, L] = 2L, [H, \Lambda] = -2\Lambda, \cr
&[d, d^*] = dd^* + d^*d = \triangle, &(6.10)\cr
&[d_c, d_c^*] = d_cd_c^* + d_c^*d_c = \triangle,\cr
&[H, d] = d, [H, d^*] = -d^*,\cr
&[H, d_c] = d_c, [H, d_c^*] = - d_c^*,\cr
&[L, d^*] = - d_c, [L, d_c^*] =  d,\cr
&[\Lambda, d] =  d_c^*, [\Lambda, d_c] =  - d^*.\cr
}$$

{\sl Theorem 6.1.}-
Let  $\g$ be
a complex finite-dimensional Lie algebra with
a non-degenerate invariant symmetric bilinear form.
Then there exist operators on
$W^{{\infty \over 2} + *}_{rel}(\tilde{\g})$, 
which are analogous to the classical operators in  K\"ahler geometry.

{\it Proof.}- 
It was shown in [FGZ] that a relative semi-infinite 
complex
$C_{\infty}^*(\l, \l_0, V)$, where $\l = \oplus_{n\in\Z}\l_n$ is a complex $\Z$-graded Lie algebra,
 and $V$
is a graded Hermitian $\l$-module, has a structure, which is similar to 
that of the de Rham complex in  K\"ahler geometry.
It is assumed that 
there exists a 2-cocycle $\gamma$ on $\l$ such that
$\gamma|_{\l_n\times \l_{-n}}$ is non-degenerate if 
$n \in \Z \backslash 0$ and it is zero otherwise.
Then there exist operators on $C_{\infty}^*(\l, \l_0, V)$ analogous to 
the classical ones.

We will define analogues of the classical operators on 
$W^{{\infty \over 2} + *}_{rel}(\tilde{\g})$. Using the form
$\langle\cdot, \cdot\rangle$  on $\g$ we obtain the 2-cocycle $\gamma$
on $\tilde{\g}$:
$$\gamma (g_1\otimes t^n, g_2\otimes t^m) =
 n\langle g_1, g_2\rangle\delta_{n, -m},
\hbox{ for } 
g_1, g_2 \in \g.\eqno (6.11)$$
Notice that
$\gamma|_{\tilde{\g}_n\times \tilde{\g}_{-n}}$ is non-degenerate if
$n \in \Z \backslash 0$ and zero otherwise.
Let
$$\Lambda^{{\infty \over 2} + *}_{rel}(\tilde{\g}) =
 \oplus_{a,b\geq 0}
\Lambda^a(\n_+')\wedge \Lambda^b_{\infty }(\n_-'). \eqno(6.12)$$
For a   homogeneous element in 
$\Lambda^a(\n_+')\wedge \Lambda^b_{\infty }(\n_-')$,
$a$ is the number of added elements, and $b$ is the number of missing elements
with respect to the vacuum  vector $\1_{rel}$.
Let
$$C^{a,b}(\tilde{\g}) = [S^{{\infty \over 2} + *}(\tilde{\g})\otimes
\Lambda^a(\n_+')\wedge \Lambda^b_{\infty }(\n_-')]^{\g_0}.\eqno (6.13)$$
We obtain a bigrading on the relative semi-infinite Weil complex, such that
$$W^{{\infty \over 2} + i}_{rel}(\tilde{\g}) = 
\oplus_{a-b=i}C^{a,b}(\tilde{\g}). \eqno (6.14)$$
Let $d$ be the restriction of the differential
 to the relative subcomplex. Notice that
$$d: C^{a,b}(\tilde{\g}) \longrightarrow
C^{a+1,b}(\tilde{\g}) \oplus  C^{a,b-1}(\tilde{\g}).\eqno (6.15)$$
Define $d_1$ and $d_2$ such that
$$\eqalignno{
&d = d_1 + d_2, &(6.16)\cr
&d_1: C^{a,b}(\tilde{\g}) \longrightarrow C^{a+1,b}(\tilde{\g}), \cr
&d_2: C^{a,b}(\tilde{\g}) \longrightarrow C^{a,b-1}(\tilde{\g}). \cr
}$$
Let
$$d_c = i(d_1 - d_2).\eqno (6.17)$$
To define the adjoint operators, we have to introduce a Hermitian form 
on $W^{{\infty \over 2} + *}_{rel}(\tilde{\g})$.

It was shown in [FGZ] that if a $\Z$-graded Lie algebra $\l$ admits an antilinear automorphism 
$\sigma$ of order 2 such that $\sigma(\l_n) = \l_{-n}$, then
there exists a 
Hermitian form on
$\Lambda^{{\infty \over 2} + *}(\l)$
such that 
$$\varepsilon (x')^* = -\varepsilon (\sigma( x')), \quad
\tau (x)^* = -\tau (\sigma (x)), 
\eqno (6.18)$$
where $x\in \l, x'\in \l'$.

To define  a Hermitian form $\lbrace \cdot, \cdot \rbrace$ on
$\Lambda^{{\infty \over 2} + *}_{rel}(\tilde{\g})$,
we set $\lbrace \1_{rel}, \1_{rel} \rbrace = 1$.
We fix a basis 
$\lbrace v_i\rbrace$  of $\g$
so that
$\langle v_i, v_j\rangle = \delta_{i,j}$. Let $u$ run through this basis.
We define an antilinear automorphism $\sigma$ of $\tilde{\g}$ as follows:
$$\sigma(u_n) = iu_{-n}. \eqno (6.19)$$
Correspondingly,
$$\sigma(u_n') = -iu_{-n}'. \eqno (6.20)$$
We introduce
a Hermitian form on $\Lambda^{{\infty \over 2} + *}_{rel}(\tilde{\g})$
so that the relations (6.18), where 
$$x\in \tilde{\g}_n, x'\in \tilde{\g}'_n \hbox{ for }n\not=0\eqno(6.21)$$
hold.
In the similar way we introduce a Hermitian form on 
$S^{{\infty \over 2} + *}(\tilde{\g})$, such that
$$\gamma (x')^* = \gamma (\sigma( x')),\quad
\beta (x)^* = -\beta (\sigma (x)). 
\eqno (6.22)$$
Then we obtain a Hermitian form $\lbrace \cdot, \cdot \rbrace$ on
$W^{{\infty \over 2} + *}_{rel}(\tilde{\g})$ by tensoring these two forms.
It gives a paring: 
$C^{a,b}(\tilde{\g}) \longrightarrow C^{b,a}(\tilde{\g})$.
 To define  a Hermitian form on 
$C^{a,b}(\tilde{\g})$, we use the linear map
$$* : C^{a,b}(\tilde{\g}) \longrightarrow C^{b,a}(\tilde{\g}),\eqno (6.23)$$
defined as follows:
$$\eqalignno{
&* \Bigl(v\otimes 
\bigl(\varepsilon (u'_{{n_1}})\cdots\varepsilon (u'_{{n_a}})
\tau (u_{{m_1}})\cdots\tau (u_{{m_b}})\1_{rel}\bigr)\Bigr) = \cr
&v\otimes\bigl(\varepsilon (u'_{{-m_1}})\cdots\varepsilon (u'_{{-m_b}})
\tau (u_{{-n_1}})\cdots\tau (u_{{-n_a}})\1_{rel}\bigr), &(6.24)\cr
}$$
where $v\in S^{{\infty \over 2} + *}(\tilde{\g})$, 
$\lbrace n_i\rbrace_{i=1}^a > 0$ and
$\lbrace m_i\rbrace_{i=1}^b < 0$.
Finally, the  Hermitian form on  $C^{a,b}(\tilde{\g})$ is defined by
$(w_1, w_2) = \lbrace i^{a+b}*w_1, w_2 \rbrace$ (cf. [FGZ]).
We introduce the adjoint operators $d^*, d_c^*$ and  the Laplace operator
$\triangle = dd^* + d^*d$.

It was pointed out in [FGZ] that as in the classical theory (see [GH]),
there exists an action of $\s\l(2)$  on  
$H^*_{\infty}(\l, \l_0, V)$. One can identify
$\l_n'$ with $\l_{-n}$ by means of the cocycle $\gamma$.
If $\lbrace e_i\rbrace$ is a homogeneous basis in $\l$,
then
$\s\l(2) = \langle L, H, \Lambda\rangle$  is defined as follows:
$$\eqalignno{ 
&L = (i/2)\sum_{m\in{\Z\backslash 0}} \varepsilon (e_m)\varepsilon (e_m'), 
&(6.25)\cr
&H = - \sum_{m\in{\Z\backslash 0}} : \tau (e_m) \varepsilon (e_m') :, \cr
&\Lambda = (i/2)\sum_{m\in{\Z\backslash 0}} \tau (e_m)\tau (e_m'). \cr
}$$
We 
identify
$\tilde{\g}_n'$ with $\tilde{\g}_{-n}$ by means of the cocycle $\gamma$
(see (6.11)),
and set
$$\E1 = L, \H1 = H,\F1 = \Lambda.\eqno (6.26)$$
Then we obtain  the 
$\S\L(2) = \langle \E1, \H1, \F1 \rangle$ defined in (5.10).
The operators
$$\lbrace\triangle, \E1, \H1, \F1,
d, d^*, d_c, d_c^*\rbrace\eqno (6.27)$$ 
are the analogues of the classical operators (6.9).
$$\eqno \SQ$$

{\sl Theorem 6.2.}-
Let  $\g$ be
a complex finite-dimensional Lie algebra with
a non-degenerate invariant symmetric bilinear form.
Then $H^{{\infty \over 2} + *}(\tilde{\g},  \tilde{\g}_0,
S^{{\infty \over 2} + *}(\tilde{\g}))$ is a module over
\hfil\break
$\hat{D}er S'(2,0)$ with central charge $3dim \g$.
 
{\it Proof.}- By Theorem 5.3, $W^{{\infty \over 2} + *}_{rel}(\tilde{\g})$
is a module over 
$\hat{D}er S'(2,0)$ with central charge $3dim \g$.
By Corollary 5.1, there is an action of
$\hat{S}'(2, 0)$  on
$H^{{\infty \over 2} + *}
(\tilde{\g}, \tilde{\g}_0,   S^{{\infty \over 2} + *}(\tilde{\g}))$.
We have proved  that
$$Der_{ext}S'(2,0) = \S\L (2) = \langle\E1, \H1, \F1\rangle, \eqno (6.28)$$
see (5.10).
Notice that as in the classical case, the element $\F1$ and the differential $d$
do not commute.
Nevertheless, there exists an action of $\S\L(2)$ on the relative semi-infinite 
cohomology according to [FGZ].
$$\eqno \SQ$$

{\sl Theorem 6.3.}- The degree zero part of the $\Z$-grading $deg$
of $S'(2,0)$ is isomorphic to the Lie
superalgebra  of classical operators in  K\"ahler geometry.

{\it Proof.}-
Recall that
the $\Z$-grading $deg$ on $S'(2,0)$ 
is defined by 
the element $L_0\in Witt$, see (4.17)-(4.19).
One can easily check that
$$S'(2,0)_{deg=0} = \langle L_0, E_{-1}, H_{0}, F_{1},
\h1_{0}^0, \p1_{0}, \x_{-1}, \y_{1}^0\rangle.\eqno (6.29)$$
The isomorphism of Lie superalgebras
$$\psi:\S \longrightarrow S'(2,0)_{deg=0} \eqno (6.30)$$
is given as follows:
$$\eqalignno{
&\psi (\triangle) = L_0,
\psi (L) = E_{-1},
\psi (H) = H_0,
\psi (\Lambda) = F_{1},&(6.31)\cr
&\psi (d) =  \h1_{0}^0,
\psi (d^*) = - \p1_{0}, 
\psi (d_c) = \x_{-1},
\psi (d_c^*) = \y_{1}^0. \cr
}$$
$$\eqno\SQ$$

{\sl Corollary 6.1.}- 
The action of $S'(2, 0)_{deg=0}$
defines a  set of quadratic
operators on
\hfil\break
$W^{{\infty \over 2} + *}_{rel}(\tilde{\g})$
(correspondingly, on $H^{{\infty \over 2} + *}(\tilde{\g},  \tilde{\g}_0,
S^{{\infty \over 2} + *}(\tilde{\g}))$),
which are analogues of the classical ones, 
and include the semi-infinite
Koszul differential $\h1 = \h1_0^0$ and the 
semi-infinite homotopy operator $\p1_0$.

{\it Remark 6.1.}-
In this work we have realized superconformal algebras by means of 
quadratic expansions on the generators of the Heisenberg and Clifford algebras
related to $\tilde{\g}$. 
Note that the differentials on a semi-infinite Weil complex
are represented by cubic expansions.
One can possibly define
an additional (to the already known) action
of the $N = 2$  SCA
on $W^{{\infty \over 2} + *}(\tilde{\g})$, considering
Fourier components of the differentials $d$ and $d^*$, [Fe].

\vskip 0.2in
   
{\sl Acknowledgements.}- This work has been partly done at the 
Max-Planck-Institut f\"ur  Mathematik in Bonn,  
L'Institut des Hautes \'Etudes Scientifiques in Bures-sur-Yvette, 
and the Institute for Advanced Study in Princeton. 
I wish to thank   MPI, IHES, and IAS 
for their hospitality and  support.
I am grateful to  B. Feigin, A. Givental, M. Kontsevich, V. Serganova,  
and V. Schechtman for very useful discussions.

\vskip 0.2in
\font\red=cmbsy10
\def\~{\hbox{\red\char'0016}}
\centerline{BIBLIOGRAPHY}
\vskip 0.2in

\item{[Ad]} M. Ademollo, L. Brink, A. D'Adda, R. D'Auria, E. Napolitano,
S. Sciuto, E. Del Giudice, P. Di Vecchia, S. Ferrara, F. Gliozzi, R. Musto and
R. Rettorino,
{\it{Dual strings with $U(1)$ colour symmetry}}, 
Nucl. Phys. B111 (1976),  77-110.

\item{[Ak]} F. Akman,
 {\it{Some cohomology operators in $2$-D field theory}},
Proceedings of the conference on Quantum topology (Manhattan, KS, 1993),
World Sci. Publ, River Edge, NJ (1994), 1-19.

\item{[Fe]} B. L. Feigin,
{\it{Private communication.}}

\item{[Fu]} D. B. Fuks,  {\it{Cohomology of infinite-dimensional
Lie algebras}}, Consultants Bureau, New York and London, 1986.

\item{[FF]} B. Feigin, E. Frenkel,  {\it{Semi-infinite Weil
Complex and the Virasoro Algebra}},
Commun. Math. Phys. 137 (1991),  617-639. Erratum: Commun. Math.
Phys. 147 (1992),  647-648.

\item{[FGZ]} I. Frenkel, H. Garland, G. Zuckerman,
 {\it{Semi-infinite cohomology and string theory,}}
 Proc. Natl. Acad. Sci. U.S.A.  83 (1986),  8442-8446.

\item{[FST]} B. L. Feigin, A. M.  Semikhatov, I. Yu. Tipunin,
{\it{Equivalence between chain categories of representations of affine
$\s\l (2)$ and $N = 2$ superconformal algebras}}, 
J. Math. Phys. 39 (1998), no 7, 3865-3905.

\item{[G]}
E. Getzler,
{\it{Two-dimensional topological gravity and equivariant cohomology,}}
Commun. Math. Phys. 163 (1994) no 3, 473-489.

\item{[GH]}
P. Griffiths, J. Harris,
{\it{Principles of algebraic geometry,}}
Wiley-Interscience Publ., New York, 1978.

\item{[KL]}
V. G. Kac, J. W. van de Leur,
{\it{On Classification of Superconformal Algebras}}, in
S. J. Gates et al., editors, {\it Strings-88},
World Scientific 1989,  77-106.

\item{[P1]}
E. Poletaeva, {\it{Semi-infinite Weil complex
 and $N = 2$ superconformal algebra I}}, preprint MPI 97-78,
{\it{Semi-infinite Weil complex
and  superconformal algebras  II}}, preprint MPI 97-79.

\item{[P2]}
E. Poletaeva, 
{\it{Superconformal algebras and Lie superalgebras of the Hodge theory}},
preprint MPI 99-136.

\item{[P3]}
E. Poletaeva, {\it{Semi-infinite cohomology and superconformal algebras}},
Comptes Rendus de l'Acad\'emie des Sciences,
t. 326, S\'erie I (1998),  533-538.

\vskip 0.1in

Elena Poletaeva

Centre for Mathematical Sciences

Mathematics, Lund University

Box 118, S-221 00 Lund, Sweden

elena$@$maths.lth.se

\bye